\documentclass{ws-ijbc}
\usepackage{psfig}

\begin{document}

\catchline{}{}{}{}{} 

\markboth{Andreas M\"uller}{Performance of Linear Adaptive Filters driven by the Ergodic Chaotic Logistic Map}

\title{On the Performance of Linear Adaptive Filters\\
driven by the Ergodic Chaotic Logistic Map}

\author{Andreas M\"uller}

\address{Chair of Mechanics and Robotics, University Duisburg-Essen\\
Lotharstrasse 1, 47057 Duisburg, Germany\\
andreas-mueller@uni-due.de}

\maketitle
\vspace{-3ex}

\begin{history}
\received{(to be inserted by publisher)}
\end{history}

\begin{abstract}
Chaotic dynamical systems are increasingly considered for use in coding and
transmission systems. This stems from their parameter sensitivity and
spectral characteristics. The latter are relevant for channel estimation
methods. In particular the logistic map  $f_\lambda =\lambda x\left(
1-x\right) $ has been employed in chaotic coding and spread spectrum
transmission systems. For $\lambda =4$ the statistical properties of
sequences generated by $f_4$ are considered as ideal drive signals for
channel estimation schemes. This assumption is proven in the present paper.
To this end the higher order statistical moments and the autocorrelation of
time series generated by $f_4$ are derived. It is shown that for $\lambda =4$
the zero mean time series is uncorrelated. The adaptation performance of
finite impulse response (FIR) digital adaptive filters (DAF) used for
channel estimation is analyzed. It is shown that using zero mean sequences
of $f_4$ leads to the maximal possible FIR DAF performance. An optimal value
for the damping parameter in the LMS scheme is derived that leads to the
maximal performance and ensures stability. The analytic considerations are
confirmed by simulation results.
\end{abstract}


\section{Introduction: Communication and Channel Estimation with Chaotic
Signals}

Chaotic coding and communication have matured since initially proposed in
the late 1980's. Chaotic communication was first proposed in \cite
{PecoraCarroll}, \cite{Carroll1995} and \cite{Parlitz1992} who showed
synchronizability of two identical chaotic systems by transmitting part of
its state, and the chaotic Chua circuit \cite{Matsumoto} was one of the
first implemented chaotic circuits. Since then different transmission
schemes were proposed employing chaotic dynamical systems for source and
channel coding, and chaotic circuits have been implemented \cite{Ditto2009},%
\cite{Senani1998}. A good overview of this topic can be found in \cite
{LauTse} and \cite{LeungYuMurali2002},\cite{tse},\cite{Yang1999}. There are
essentially two attributes that make chaotic systems attractive for
communication purposes: the inherent parameter sensitivity and the
statistical properties of the generated chaotic signals. The first is
relevant for source coding and chaotic cryptography \cite{Alvarez},\cite
{Kocarev},\cite{Kocarev2001},\cite{Martinez-GuerraYu2008},\cite
{MiLiaoChen2007},\cite{ParlitzKocarev1996},\cite{WangWangPei2011}, while the
latter is relevant for channel coding with potential for spread spectrum
communication (aiming on maximizing the bandwidth allocation) \cite
{HaiJiandong1997},\cite{Mazzini},\cite{Mazzini2},\cite{Tou},\cite
{YangChua2000} and for channel estimation (monitoring of the transmission
channel) \cite{ICC1999},\cite{Feng2004}.

Spread spectrum communication and channel estimation methods take advantage
of the statistical properties of the transmitted signal rather than its
deterministic chaotic nature. Spread spectrum communication relies on
so-called spreading codes, and chaotic spreading codes could offer
complementary alternatives \cite{Vitali},\cite{Ye}. Channel estimation
strategies commonly infer the transmission characteristics in a non-blind
fashion, i.e. by comparing the original signal with that received at the
channel output. Now the performance of channel estimation algorithms depends
on the higher order statistical moments of the monitored signal, and is
maximal if the signal exhibits a white spectrum (flat power spectral
density) \cite{haykin2001}. In practice channel estimation algorithms
operate on-line, i.e. use the transmitted signal when the channel is in
operation. Due to the generally inappropriate spectral properties of the
transmitted information signal the update performance is low. Therefore
off-line adaptation is frequently pursued by injecting tailored test
signals. This apparently intrudes the actual transmission process, however.
Since chaotic sequences possess advantageous statistical properties chaotic
channel coding (encoding the signal to be transmitted in a chaotic carrier
sequence) potentially allows for on-line channel estimation with increased
update performance. Aiming on maximal update performance it is crucial to
analyze the statistical properties of the chaotic sequence. Different
chaotic systems may have drastically varying spectral properties, and it is
well-known that certain chaotic systems have an almost white spectrum while
others (like the Lorenz system) exhibit a colored spectrum.

The one-dimensional logistic map (LM) 
\begin{equation}
f_\lambda \left( x\right) =\lambda x\left( 1-x\right)  \label{logisticmap}
\end{equation}
exhibits chaotic dynamics for parameter values $\lambda =3.57,\ldots ,4$,
and has been the subject of extensive research. A variety of chaotic coding
schemes were development where the corresponding time sequence $\{x_i\}$ is
obtained by iterating the difference equation $x_i=f_\lambda \left(
x_{i-1}\right) $, and the actual value of $\lambda $ is related to the
information to be transmitted. In this fashion the LM has been used for
chaotic encryption \cite{Alvarez},\cite{chen},\cite{Kocarev},\cite
{Kocarev2001},\cite{MiLiaoChen2007},\cite{PatidarSu2009}, for chaotic
modulation and masking \cite{ParlitzKocarev1996},\cite{t-com2002}, as
carrier signal for wide-band communication \cite{KolumbanKrebesz2009},\cite
{HaiJiandong1997}, and for channel estimation \cite{ICC1999}.

It was conjectured that the statistical properties of LM are ideal for
channel estimation if the parameter value is close to $\lambda =4$. However,
this has so far only been verified by numerical simulations, e.g. \cite
{PatidarSu2009}.

In this paper the adaptation performance of finite impulse response (FIR)
digital adaptive filters (DAF) is analyzed when driven by chaotic time
series generated by the ergodic LM with $\lambda =4$. To this end the higher
order statistical moments and the autocorrelation of such time series are
derived analytically in section \ref{secLM}. It is shown in particular that
the time series for $\lambda =4$ is uncorrelated and exhibits a flat power
spectral density. In order to establish the FIR DAF performance the average
coefficient vector and its fluctuation from this average is analyzed in
section 3 when the standard least mean square (LMS) algorithm is applied for
adaptation. It is concluded that the maximum adaptation performance of FIR
DAF is achieved for zero-mean sequences generated by $f_4$, which approaches
the maximal possible performance of white Gausian noise. Numerical examples
are presented confirming the theoretical result for $f_4$, and showing the
performance attenuation experienced when applying the LM with $\lambda <4$.

\section{%
\label{secLM}%
Characteristics of the Ergodic Chaotic Logistic Map $f_4$}

\subsection{Invariant Density}

\begin{figure}[tbh]
\centerline{\psfig{file=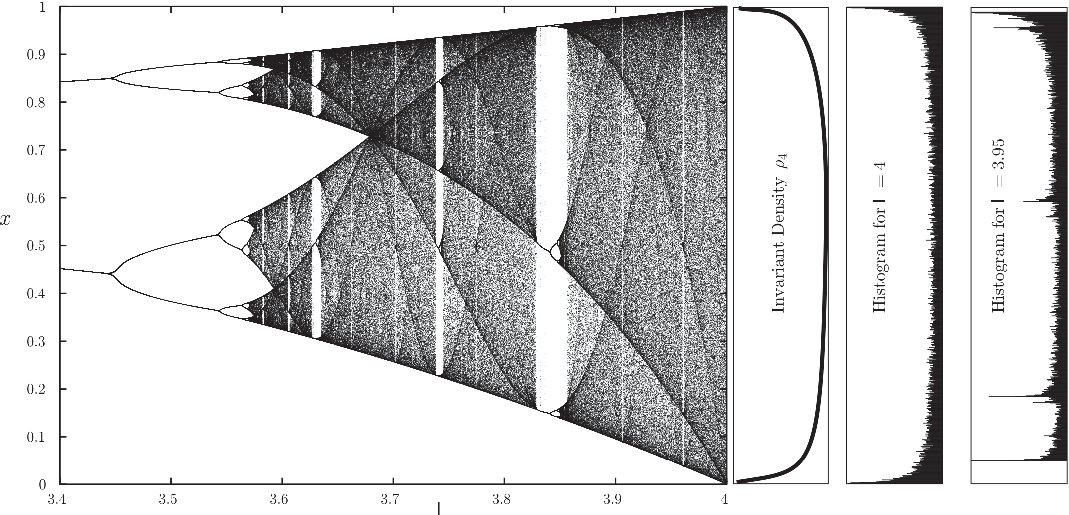,width=18cm 
} }
\caption{Bifurcation diagram of the Logistic Map, the invariant density $%
\rho _4$ for the ergodic $f_4$, and the histograms for $10^5$ samples of the
LM with $\lambda =4$ and $\lambda =3.95$.}
\label{log_dens}
\end{figure}

An important aspect for chaotic coding is the distribution of the sequence $%
\{x_i\}$. Moreover, the qualitative study of a discrete dynamical system is
based on the concept of invariant sets. $U\subset {\Bbb{R}}$ is an invariant
set of $x_i=f_\lambda \left( x_{i-1}\right) $ if $x_i\in U$ for any $k$ and $%
x_0\in U$. Figure \ref{log_dens} shows the bifurcation diagram of the LM for 
$\lambda =3.4,\ldots ,4$. For $\lambda =4$ the invariant set is the whole
interval $\left[ 0,1\right] $, whereas for $\lambda <4$ the invariant sets
are contained in $\left[ 0,1\right] $. Information about the distribution
within an invariant set is given by the invariant density that determines
the density of iterates when started at some point $x_0$. If this density
does not depend on the starting point $x_0$, the system is called ergodic,
and the invariant density is denoted $\rho _\lambda \left( x\right) $. It
was shown in \cite{ulam},\cite{ott} that the LM becomes ergodic for $\lambda
=4$ and that 
\begin{equation}
\rho _4\left( x\right) =\frac 1\pi \frac 1{\sqrt{x\left( 1-x\right) }}.
\label{rho}
\end{equation}
The invariant density is the equivalence of the invariant probability
distribution of a stochastic process. Ergodicity implies that orbits of $f_4$
generated from arbitrary initial points approach any point in the strange
attractor after an infinite number of iterations. This has important
implications for the coding and signal transmission, and it is desirable for
various reasons that the sequence $\{x_i\}$ is symmetrically distributed at
least. For instance from a technical point of view it is important that the
sequence being transmitted can be converted to a zero-mean signal, which is
possible using the ergodic LM (since it has the constant mean value 0.5).
Also for instance the cryptosystem proposed in \cite{Baptista1998} relies on
the ergodicity of the employed chaotic dynamical system. The LM being
ergodic implies that any invariant subset of $\left[ 0,1\right] $ has
measure zero with respect to $\rho _4$. One also says that $f_4$ is $\rho _4$%
-ergodic.

The ergodicity ceases for $\lambda <4$. Moreover the distribution and its
symmetry changes drastically for $\lambda <4$. The asymmetry is apparent
from the histogram for $\lambda =3.95$ in figure \ref{log_dens}. Eventually
for several intervals of $\lambda $ periodic windows occur where the system
is not chaotic. Furthermore the variation of bifurcation parameter $\lambda $
has consequences for the statistical properties at large.

The update performance of adaptive filters is determined by the statistical
moments of the drive signal, be it generated by a random or deterministic
process. Therefore, as a prerequisite, the analytic expression (\ref{rho})
will be used to derive the higher order statistical moments and the
autocorrelation of sequences generated by $f_4$.

\newpage%

\subsection{Statistical Moments}

The statistical moments $E\left[ x_t^\nu \right] $ of a stochastic process $%
x_t$ ($E$ denotes the expectation) can be expressed in terms of the
characteristic function (the inverse Fourier transform) of its probability
density function. This can be adopted for ergodic chaotic systems by
replacing the probability density by the invariant density. Then the $\nu $%
th order statistical moment is given by 
\begin{equation}
m_\nu :=E\left[ x_t^\nu \right] =\lim_{\xi \rightarrow 0}\frac 1{i^\nu }%
\frac{d^\nu }{d\xi ^\nu }F_4\left( \xi \right) ,  \label{moment1}
\end{equation}
where $i$ is the imaginary unit and 
\begin{equation}
F_4\left( \xi \right) =\frac 1\pi \int_0^1\frac{e^{ix\xi }}{\sqrt{x\left(
1-x\right) }}dx  \label{charfunct}
\end{equation}
is the characteristic function of $\rho _4\left( x\right) $. The expression (%
\ref{moment1}) can be evaluated noticing the similarity between $F_4$ and
the Kummer function of the first kind defined as \cite{Bronstein} 
\begin{equation}
\Psi \left( a,b,\xi \right) =\frac{\Gamma \left( c\right) }{\Gamma \left(
a\right) \Gamma \left( c-a\right) }\int_0^1x^{a-1}\left( 1-x\right)
^{c-a-1}e^{\xi x}dx  \label{kummer1}
\end{equation}
for $b>a>0$. The Kummer function admits the series expansion 
\begin{equation}
\Psi \left( a,b,\xi \right) =\sum_{r=0}^\infty \frac{\left( a\right) _r}{%
\left( b\right) _rr!}\xi ^r  \label{kummer2}
\end{equation}
where 
\begin{equation}
\left( a\right) _0=1,\;\left( a\right) _r=a\left( a+1\right) \left(
a+2\right) \cdots \left( a+r-1\right) .
\end{equation}
Since $\Psi \left( a,b,\xi \right) $ is analytic in $\Bbb{R}$ and converges
for all $\xi $ it may be analytically continued in the complex plane $\Bbb{C}
$. Hence, noting (\ref{kummer1}), the characteristic function (\ref
{charfunct}) can be rewritten as 
\begin{equation}
F_4\left( \xi \right) =\frac 1\pi \frac{\Gamma \left( a\right) \Gamma \left(
c-a\right) }{\Gamma \left( c\right) }\Psi \left( a,b,i\xi \right) 
\end{equation}
with $a=\frac 12$ and $b=1$, and $i$ the imaginary unit. It can be verified
from (\ref{kummer2}) that the $\nu $th derivation of $\Psi $ satisfies 
\begin{equation}
\frac{d^\nu }{d\xi ^\nu }\Psi \left( a,b,\xi \right) =\frac{\left( a\right)
_\nu }{\left( b\right) _\nu }\Psi \left( a+\nu ,b+\nu ,\xi \right) .
\label{kummerdif}
\end{equation}
Combining (\ref{kummerdif}) and (\ref{kummer2}) leads to the following
expression for the statistical moments 
\begin{eqnarray}
m_\nu  &=&\lim_{\xi \rightarrow 0}\frac 1\pi \frac{\left( a\right) _\nu }{%
\left( b\right) _\nu }\frac{\Gamma \left( a\right) \Gamma \left( c-a\right) 
}{\Gamma \left( c\right) }\Psi \left( a+\nu ,b+\nu ,\xi \right)   \nonumber
\\
&=&\frac 1\pi \frac{\left( a\right) _\nu }{\left( b\right) _\nu }\frac{%
\Gamma \left( a\right) \Gamma \left( c-a\right) }{\Gamma \left( c\right) }%
\frac{\left( a+\nu \right) _0}{\left( b+\nu \right) _0}=\frac 1\pi \frac{%
\left( a\right) _\nu }{\left( b\right) _\nu }\frac{\Gamma \left( a\right)
\Gamma \left( c-a\right) }{\Gamma \left( c\right) }.
\end{eqnarray}
Finally inserting $a=\frac 12$ and $c=1$ yields the $\nu $th order moment of
sequences generated by the LM $f_4$ 
\begin{equation}
m_\nu =\frac{\left( \frac 12\right) _\nu }{\nu !}\frac{\Gamma ^2\left( \frac
12\right) }\pi =\frac{1\cdot 3\cdot 5\cdot \ldots \cdot \left( 2\nu
-1\right) }{2^\nu \nu !}.  \label{mn}
\end{equation}
The first 7 moments are listed in table 1. Clearly, $m_0=1$ means that
certainly a $x_i$ takes place in $\left[ 0,1\right] $ whereas $m_1=\frac 12$
is the expectation of the sequence $\left\{ x_i\right\} $.

\begin{table}[htbp] \centering%
\tbl{Statistical moments of the logistic map $f_4$} {\ 
\begin{tabular}{|l||l|l|l|l|l|l|l|l|}
\hline
order $\nu $ & $0$ & $1$ & $2$ & $3$ & $4$ & $5$ & $6$ & $7$ \\ \hline
$m_\nu $ & $1$ & $\frac 12$ & $\frac 38$ & $\frac 5{16}$ & $\frac{35}{128}$
& $\frac{63}{256}$ & $\frac{231}{1024}$ & $\frac{429}{2048}$ \\ \hline
\end{tabular}
} \label{momtab1} 
\end{table}%

In technical implementations it is desirable to use a zero-mean chaotic
sequence (electrically a DC free signal), obtained by subtraction of the
mean value $m_1$. Moreover, it will become clear in section \ref{FIR} that
zero-mean sequences lead to the best performance of channel estimation
schemes. Subtracting the mean value $\overline{x}=m_1=\frac 12$, and using $%
E[\left( x-\overline{x}\right) ^2]=E\left[ x^2\right] -E\left[ x\right] ^2$
and so forth yields the $\nu $th order moments $\widetilde{m}_\nu $ for the
zero-mean sequence in table 2. Unlike for the moments $m_\nu $ there is no
obvious systematic rule for expressing $\widetilde{m}_\nu $ such as (\ref{mn}%
).

\begin{table}[htbp] \centering%
\tbl{Moments of the zero-mean sequences of the logistic map ${f}_4$} {\ 
\begin{tabular}{|l||c|c|c|c|c|c|c|c|c|c|c|c|c|c|c|}
\hline
order $\nu $ & $0$ & $1$ & $2$ & $3$ & $4$ & $5$ & $6$ & $7$ & $8$ & $9$ & $%
10$ & $11$ & $12$ & $13$ & $14$ \\ \hline
$\widetilde{m}_\nu $ & $1$ & $0$ & $\frac 18$ & $0$ & $\frac 3{128}$ & $0$ & 
$\frac 5{1024}$ & $0$ & $\frac{35}{32768}$ & $0$ & $\frac{63}{262144}$ & $0$
& $\frac{231}{4194304}$ & $0$ & $\frac{429}{33554432}$ \\ \hline
\end{tabular}
} \label{momtab2} 
\end{table}%

\subsection{Autocorrelation}

The performance of digital channel estimation algorithms depends on the
eigenvalue spread of the autocorrelation matrix of the driving signal \cite
{haykin2001}, and the best performance is obtained for uncorrelated samples.
The correlation of sequences generated by chaotic systems is also
significant for chaotic coding \cite{Deane2000}. In this section the
autocorrelation function will be derived for the LM $f_4$. In this paper the
autocorrelation matrix of a process $\{x_i\}$ is used as defined in signal
processing: $R_{ij}=E\left[ x_ix_j\right] $ since this appears in the
analysis in section \ref{FIR}. This only differs from the definition in
statistics in that the samples are normalized by subtraction of the mean
value.

The autocorrelation indicates the relationship of sample values $x_i$ and $%
x_{i+m}$. The samples originate from an initial state $x_0$ by application
of $f_4$, $i$ times and $i+m$ times, respectively, so that 
\begin{equation}
x_{i+m}=f_4^{i+m}\left( x_0\right) =f_4^m\circ f_4^i\left( x_0\right)
=f_4^m\left( x_i\right)
\end{equation}
Since $f_4$ is ergodic the correlation only depends on the distance $m$, and
the invariant measure $\rho _4$ exists that gives rise to the
autocorrelation defined as (see also \cite{schuster2005}) 
\begin{equation}
C_4\left( m\right) =\int_0^1\rho _4\left( x\right) xf_4^m\left( x\right)
dx=\frac 1\pi \int_0^1\frac{f_4^m\left( x\right) x}{\sqrt{x\left( 1-x\right) 
}}dx,  \label{cm}
\end{equation}
which is formally equivalent to that of a stationary stochastic process when
the probability distribution is replaced by the invariant density. The
integral can be solved in closed from as 
\begin{equation}
\frac 1\pi \int_0^1\frac{f_4^m\left( x\right) x}{\sqrt{x\left( 1-x\right) }}%
dx=\sum_{k=0}^{m+1}\frac{\alpha _k}\pi x^k\sqrt{x-x^2}+\left. \frac{\beta _k}%
\pi \arcsin \left( 2x-1\right) \right| _{x=0}^{x=1}  \label{int}
\end{equation}
with certain $\alpha _k\in {\Bbb{R}}$ and 
\[
\beta _k=\left\{ 
\begin{array}{ll}
\frac 38,\;\; & m=0%
\vspace{-1ex}%
\\ 
&  \\ 
\frac 14, & m\geq 1\;\;.
\end{array}
\right. 
\]
The first term on the right hand side of (\ref{int}) clearly vanishes for $%
x=0,1$. Thus the autocorrelation (\ref{cm}) is explicitly 
\begin{equation}
C_4\left( m\right) =\frac 14+\frac 18\delta _{m,0}  \label{C4}
\end{equation}
where $\delta _{i,j}$ is the Kronecker symbol. In the later application the
zero-mean chaotic sequence $\left\{ x_i-\overline{x}\right\} $ will be used.
For this sequence the autocorrelation function is 
\begin{eqnarray}
\widetilde{C}_4\left( m\right) &=&\int_0^1\rho _4\left( x\right) \left( x-%
\overline{x}\right) \left( f_4^m\left( x\right) -\overline{x}\right)
dx=\int_0^1\rho _4\left( x\right) xf_4^m\left( x\right) dx-\overline{x}^2 
\nonumber \\
&=&C_4\left( m\right) -\frac 14
\end{eqnarray}
which together with (\ref{C4}) yields 
\begin{equation}
\widetilde{C}_4\left( m\right) =\frac 18\delta _{m,0}.  \label{C40}
\end{equation}

Hence samples of the zero-mean chaotic sequence $\left\{ x_i-\overline{x}%
\right\} $ are uncorrelated, whereas the samples of $\left\{ x_i\right\} $
are correlated. Consequently only the zero-mean sequence exhibits a flat
power spectrum.

\begin{remark}
It must be noticed that these results strictly only apply to $f_4$.
Numerical simulations show that if the bifurcation parameter $\lambda $ is
varied, the signal becomes correlated and the spectral density is no longer
flat. This is one reason why within chaotic coding schemes $\lambda $ is
kept close to $4$. That is, simultaneous variation of $\lambda $ and
flattening the power spectrum is only possible in a very narrow $\lambda $%
-range. Further figure \ref{log_dens} shows that for decreasing values of $%
\lambda <4$ the energy distribution within the signal range becomes
unbalanced, and the orbit plot gets its first window (showing a zero
Lyaponov exponent of $f$) so that the chaotic band of $f_4\left( x_n\right) $
is split into two subbands that further split for decreasing $\lambda $
according to the Feigenbaum attractor \cite{schuster2005}. This can be
critical for technical applications. To extend the range of admissible $%
\lambda $ values, in \cite{chen} a modified LM is proposed that yields a
uniformly distributed sequence in the range $\left[ 0,1\right] $.
\end{remark}

\section{%
\label{FIR}%
Performance of FIR DAF driven by $f_4$}

\begin{figure}[tbh]
\centerline{\psfig{file=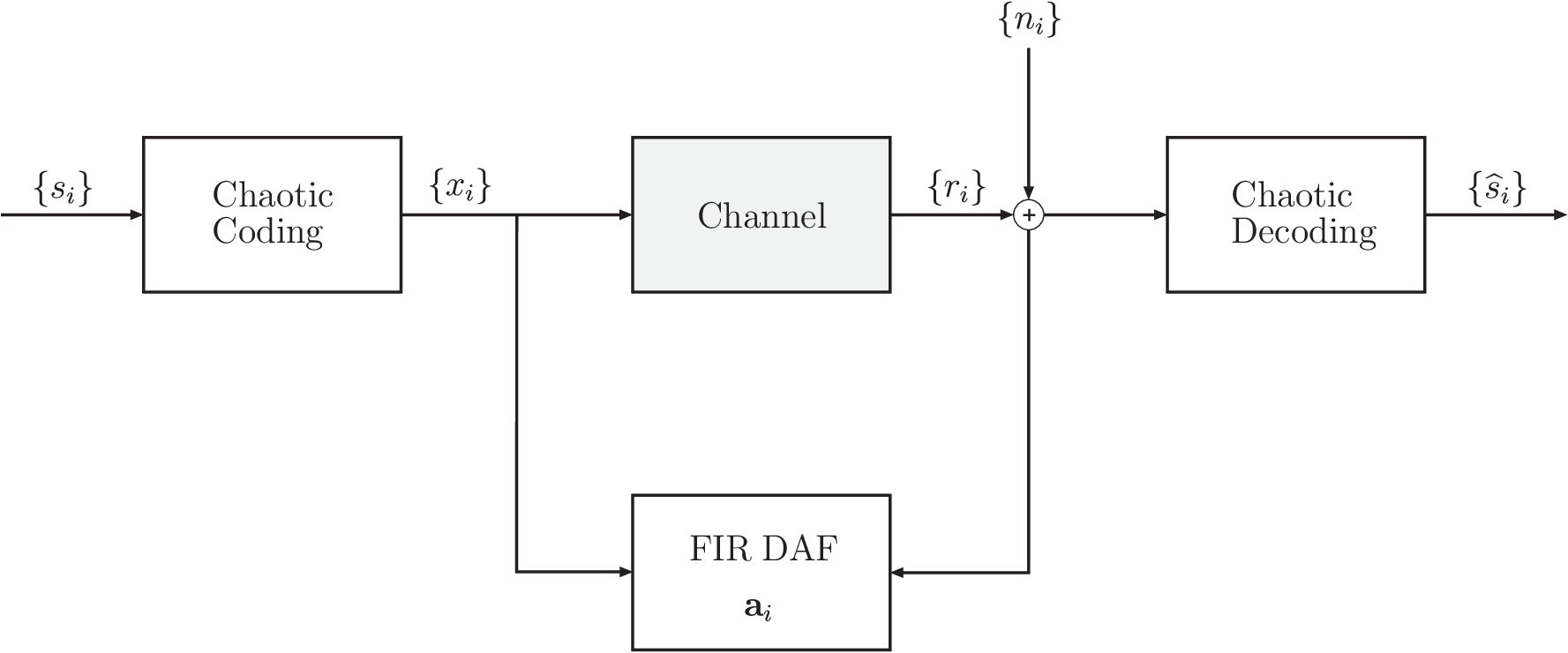,width=13.0cm}}
\caption{FIR DAF channel estimation concept using chaotic coding. The
information ${s_i}$ is encoded in the chaotic sequence ${x_i}$. The signal ${%
r_i}$ received at the channel output is decoded to reconstruct an
information signal $\widehat{s}_i$. The FIR coefficient vector $\mathbf{a}_i$%
, estimating the channel transfer function, is updated by an LMS scheme
based on the error $e_i=r_i-\mathbf{a}_i^T{x}_i$. }
\label{figFIRDAF}
\end{figure}

Channel estimation, or equalization, is an important topic for digital
communication systems and gained increasing importance with the rise of
mobile systems. A classical but still widely used approach to channel
estimation is based on FIR digital adaptive filters (DAF) \cite{goodwin},%
\cite{haykin2001}. The unknown channel is represented by a FIR DAF whose tap
weight coefficients are adapted by a least mean squares (LMS) algorithm. It
is well-known from the FIR DAF theory that the update performance depends on
the spectrum of the autocorrelation matrix and the fourth order moments of
the drive signal. If the information is encoded in a chaotic sequence
(figure \ref{figFIRDAF}), which is transmitted and employed for estimation
the transmission channel, the statistics of the chaotic sequence are
decisive. Now if the LM $f_4$ is used to generate the transmitted signal,
the results in section \ref{secLM} admit analyzing the FIR DAF performance.

In this section an upper bound on the FIR DAF performance is derived when
driven by sequences generated by $f_4$. Now chaotic coding schemes modulate
the value of $\lambda $. However, it is known from numerical simulations 
\cite{PatidarSu2009} that the samples become non-uniformly distributed in
the range $\left[ 0,1\right] $, which causes a moving mean value and
variable power spectrum. Therefore the variation of $\lambda $ is limited,
and hence results for $\lambda =4$ may be considered as a practically
sensible estimate.

\subsection{LMS Scheme for Digital Adaptive Filters}

A FIR DAF of length $m$ is parameterized by the vector $\mathbf{a}_i=\left(
a_0,\ldots ,a_m\right) $ of tap weight coefficients. Denote with $\mathbf{x}%
_i=\left( x_i,\ldots ,x_{i-m}\right) $ the vector comprising the input
values at the preceding $m$ time steps. The output of the FIR DAF model at
time step $i$ is 
\begin{equation}
\widehat{r}_i=\mathbf{a}_i^T\mathbf{x}_i.
\end{equation}
Let $\mathbf{b}$ be the vector comprising the first $m+1$ samples of the
impulse response of the unknown transmission channel (i.e. an FIR
approximation of length $m$ ). The actual channel output is $r_i=\mathbf{b}^T%
\mathbf{x}_i+\xi _i$, where $\xi _i$ represents the residuum due to the
finite approximation, and the adaptation error of the FIR DAF is $e_i=r_i-%
\widehat{r}_i=\left( \mathbf{b}-\mathbf{a}_i\right) ^T+\xi _i+n_i$, where $%
n_i$ represents additive measurement noise. The adaptation is indeed
required to converge to an optimal coefficient vector $\mathbf{a}$. Assuming
a stationary input sequence $x_i$, white Gausian measurement noise, and
independence of the terms in $e_i$ (which makes sense for slow adaptation),
the optimal tap weight coefficient vector is given by the Wiener solution 
\cite{widrow}: $\mathbf{a}_{\text{opt}}=\mathbf{b}+\mathbf{R}^{-1}\mathbf{%
\rho }_i$, with autocorrelation matrix $\mathbf{R}=E\left[ \mathbf{x}_i%
\mathbf{x}_i^T\right] $, and $\mathbf{\rho }_i=E\left[ \xi _i\mathbf{x}%
_i^T\right] $. A (LMS) method is commonly used that iteratively minimizes $%
E\left[ e_i^2\right] $ with a stochastic steepest descent update 
\begin{eqnarray}
\mathbf{a}_{i+1} &=&\mathbf{a}_i-\mu \nabla E\left[ e_i^2\right]  \nonumber
\\
&=&\left( \mathbf{I}-\mu \mathbf{R}\right) \mathbf{a}_i+\mu \left( \mathbf{Rb%
}+\mathbf{\rho }_i\right)  \label{LMS1}
\end{eqnarray}
where $\mu $ is a damping parameter. Assuming uncorrelated $\xi _i$ and $x_i$%
, the LMS update algorithm in terms of the samples $x_i$ is \cite{haykin2001}
\begin{equation}
\mathbf{a}_{i+1}=\mathbf{a}_i+\mu e_i\mathbf{x}_i\text{.}  \label{LMS2}
\end{equation}
The damping parameter determines the speed of convergence but must be
bounded by the signal power so to ensure stability: $\mu \leq 1/\left\| 
\mathbf{x}_i\right\| ^2$.

\subsection{Average Coefficient Vector}

The FIR DAF LMS scheme is introduced for general stationary drive signals,
and its performance can be analyzed taking into account the signal
statistics \cite{haykin2001}. Now the LMS performance can be analyzed
explicitly if excited by sequences generated by the chaotic LM $f_4$.

A basic performance characteristic is the average coefficient vector. It is
crucial that this average converges to the unknown impulse response $\mathbf{%
b}$, despite fluctuations due to the update process. With (\ref{LMS1}) the
expectation of the filter coefficients at step $i$ can be approximated as 
\begin{equation}
E\left[ \mathbf{a}_i\right] =\left( \mathbf{I}-\mu \mathbf{R}\right) E\left[ 
\mathbf{a}_{i-1}\right] +\mu \left( \mathbf{Rb}+\mathbf{\rho }_i\right)
\end{equation}
presuming uncorrelated $\mathbf{x}_i$ and $\mathbf{a}_i$, and stationary $%
\{x_i\}$. Subtraction of the optimal coefficient vector $\mathbf{a}_{\text{%
opt}}$ gives an error $\mathbf{\varepsilon }_i=E\left[ \mathbf{a}_i\right] -%
\mathbf{b}-\mathbf{R}^{-1}\mathbf{\rho }_i$. This error is recursively
expressed as $\mathbf{\varepsilon }_i=\left( \mathbf{I}-\mu \mathbf{R}%
\right) ^i\mathbf{\varepsilon }_0$, when the iteration starts from an
initial error $\mathbf{\varepsilon }_0$. The decay of $\mathbf{\varepsilon }%
_i$ depends solely on $\mu $ and $\mathbf{R}$, and thus on the eigenvalue
spectrum of the autocorrelation matrix, which is a positive definite
Toeplitz matrix. Denote with $\lambda _k$ the eigenvalues of $\mathbf{R}$,
the condition for $\mathbf{\varepsilon }_i$ decaying exponentially to zero
is that $-1<1-\mu \lambda _k<1$. Since $\lambda _k$ are all real and
positive, the condition for decaying $\mathbf{\varepsilon }_i$ is that 
\begin{equation}
0<\mu <2/\lambda _{\max }  \label{lambdamax}
\end{equation}
denoting with $\lambda _{\max }$ the largest eigenvalue. The error $\mathbf{%
\varepsilon }$ can be split according to the eigenvectors, and the $k$th
mode of $\mathbf{\varepsilon }$ decays proportional to $1-\mu \lambda _k$.
For $\mu =2/\lambda _{\max }$, which yields the fasted convergence, the
slowest converging mode converges proportional to $\left| 1-2\lambda _{\min
}/\lambda _{\max }\right| $. Hence, the uniformity of convergence depends on
the condition number $\sigma =\lambda _{\max }/\lambda _{\min }$ of $\mathbf{%
R}$. This is considered in the following separately for $\{x_i\}$ and $\{x_i-%
\overline{x}\}$%
\newpage%

\begin{enumerate}
\item  Original Sequence $\{x_i\}$ generated by $f_4$\newline
The $\left( m+1\right) \times \left( m+1\right) $ autocorrelation matrix for
sequences generated by the LM $f_4$ follows from (\ref{C4}) as 
\begin{equation}
\mathbf{R}=\left( 
\begin{array}{crrrrc}
3/8 & \;\;\;\;1/4 & \;\;1/4 & \;\;\;\cdots  &  & 1/4%
\vspace{1ex}%
\\ 
1/4 & 3/8 & 1/4 & \cdots  & \;\; & 1/4 \\ 
1/4 & 1/4 & \ddots  &  &  & \vdots  \\ 
\vdots  & \multicolumn{1}{c}{\vdots } & \multicolumn{1}{c}{} & 
\multicolumn{1}{c}{} & \multicolumn{1}{c}{} &  \\ 
1/4 & 1/4 & \cdots  &  &  & 
\end{array}
\right) .
\end{equation}
Its eigenvalues are $\lambda _0=\frac{2m+1}8$ and $\lambda _i=\frac
18,i=1,\ldots ,m$, so that the eigenvalue spread of $\mathbf{R}$ is $\sigma
=2m+1$. Hence the convergence of the fastest mode is $2m+1$ times faster
than that of the slowest. Furthermore the update performance of the FIR DAF
decreases linearly with the filter length $m$.

\item  Zero-Mean Sequence $\{x_i-\overline{x}\}$\newline
If a zero-mean sequence generated by $f_4$ is employed as drive signal, (\ref
{C40}) shows that the autocorrelation matrix is 
\begin{equation}
\mathbf{R}=\frac 18\mathbf{I}
\end{equation}
with eigenvalues $\lambda _i=\frac 18,i=0,\ldots ,m$. Hence the eigenvalue
spread of $\mathbf{R}$ is $\sigma =1$, which implies a uniform adaptation
performance.
\end{enumerate}

It can thus be concluded that the zero-mean chaotic sequence leads to a fast
and uniform convergence.

\subsection{Fluctuation of the Average Coefficient Vector}

Besides the convergence to the optimal Wiener solution the fluctuation of
the FIR DAF coefficient vector from the optimum, $\mathbf{b}$, is crucial,
which determines stability of the update scheme. This fluctuation can be
quantified by $\left\| \mathbf{a}_{i+1}-\mathbf{b}\right\| ^2$, and the
expectation of the fluctuation can be recursively expressed using (\ref{LMS1}%
) as 
\begin{equation}
E\left[ \left\| \mathbf{a}_{i+1}-\mathbf{b}\right\| ^2\right] =E\left[
\left( \mathbf{a}_i-\mathbf{b}\right) ^TE\left[ \left( \mathbf{I}-\mu 
\mathbf{x}_i\mathbf{x}_i^T\right) ^2\right] \left( \mathbf{a}_i-\mathbf{b}%
\right) \right] +\mu ^2E\left[ \xi _i^2\right] E\left[ \mathbf{x}_i^T\mathbf{%
x}_i\right] .  \label{exp}
\end{equation}
Assuming that $\mathbf{a}_i$, $\xi _i$, and the drive signal $x_i$ are
uncorrelated, the middle term in (\ref{exp}) becomes 
\begin{equation}
E\left[ \left( \mathbf{I}-\mu \mathbf{x}_i\mathbf{x}_i^T\right) ^2\right] =%
\mathbf{I}-2\mu E\left[ \mathbf{x}_i\mathbf{x}_i^T\right] +\mu ^2E\left[
\left( \mathbf{x}_i^T\mathbf{x}_i\right) \left( \mathbf{x}_i\mathbf{x}%
_i^T\right) \right] .  \label{mat1}
\end{equation}
This matrix determines the decay of the fluctuation. The autocorrelation
matrix $\mathbf{R}=E\left[ \mathbf{x}_i\mathbf{x}_i^T\right] $ has been
derived above, and it remains to determine $E\left[ \left( \mathbf{x}_i^T%
\mathbf{x}_i\right) \left( \mathbf{x}_i\mathbf{x}_i^T\right) \right] $. The
required 4th moments can be evaluated in closed form as it was done in (\ref
{int}). This computation is straightforward and the details are omitted.

\begin{enumerate}
\item  Original Sequence $\{x_i\}$ generated by $f_4$\newline
As for the autocorrelation the matrix, $E\left[ \left( \mathbf{x}_i^T\mathbf{%
x}_i\right) \left( \mathbf{x}_i\mathbf{x}_i^T\right) \right] $ turns out to
be dense and there is no obvious rule for the matrix elements that depend on 
$m$. Therefore, and since the sequence $\{x_i\}$ leads to non-uniform
convergence, it is not considered in further detail here. The condition for
exponential convergence of the fluctuation to zero is that the eigenvalues
of (\ref{mat1}) are less than 1. If the matrix (\ref{mat1}) were determined
for specific $m$, this condition gives rise to an upper bound for the
damping parameter $\mu $.

\item  Zero-Mean Sequence $\{x_i-\overline{x}\}$\newline
For the zero-mean sequence all off-diagonal 4th order moments vanish leading
to the diagonal $\left( m+1\right) \times \left( m+1\right) $ matrix 
\begin{equation}
E\left[ \left( \mathbf{x}_i^T\mathbf{x}_i\right) \left( \mathbf{x}_i\mathbf{x%
}_i^T\right) \right] =\left( \frac 3{128}+\frac m{64}\right) \mathbf{I.}
\label{corr1}
\end{equation}
\newline
Thus $E\left[ \left\| \mathbf{a}_{i+1}-\mathbf{b}\right\| ^2\right] =\left(
1-\frac 14\mu +\left( \frac 3{128}+\frac m{64}\right) \mu ^2\right) E\left[
\left\| \mathbf{a}_i-\mathbf{b}\right\| ^2\right] +\mu ^2E\left[ \xi
_i^2\right] $. The condition for exponential decay requires that $\left|
1-\frac 14\mu +\left( \frac 3{128}+\frac m{64}\right) \mu ^2\right| <1$,
which gives rise to a further bound on $\mu $. The minimal absolute value is
attained for $\mu _{\max }=16/\left( 3+2\,m\right) $, while $\mu _{\min }=0$
leads to the maximum for which no convergence is achieved. This is a more
stringent condition than $\mu <2/\lambda _{\max }=16$ dictated by (\ref
{lambdamax}). The final condition on the damping parameter is thus 
\begin{equation}
0<\mu \leq \frac{16}{3+2\,m}  \label{mumax}
\end{equation}
for an FIR DAF of order $m$. Clearly $\mu $ tends to zero for large $m$,
which is in accordance with the FIR DAF theory.
\end{enumerate}

\section{Simulation Results}

\begin{figure}[b]
\centerline{\psfig{file=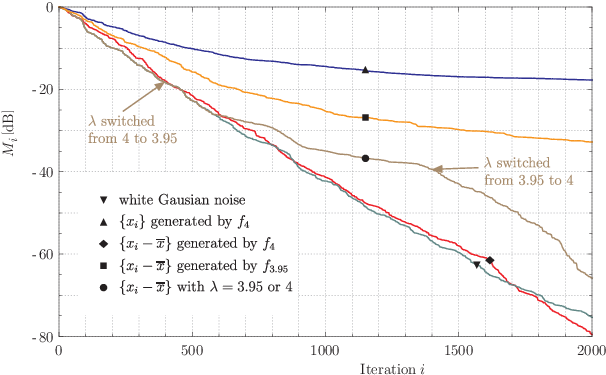,width=13.0cm}}
\caption{Adaptation performance, measured by the MMA $M_i=\left\| \mathbf{a}%
_i-\mathbf{b}\right\| ^2$ in dB, when the FIR DAF is driven by white Gausian
noise, a sequence $\left\{ x_i\right\} $ generated by the LM $f_4$, the
corresponding zero-mean sequence $\{x_i-\overline{x}\}$, and a sequence $%
\{x_i-\overline{x}\}$ generated by the LM $f_{3.95}$. Also shown is the
effect of switching $\lambda $ from 4 to 3.95 after 400 iterations and back
to 4 after 1400 iterations.}
\label{figadapt}
\end{figure}

To illustrate the analytical results the LMS update has been simulated for
an FIR DAF of length $m=128$ applied to estimate a channel modeled by an IIR
transfer function 
\[
H\left( z\right) =\frac 1{1-0.2z^{-1}+0.49z^{-2}+0.292z^{-3}}. 
\]
The FIR representation of this channel consists of the first 128 samples of
the impulse response constituting the filter coefficient vector $\mathbf{b}$.

The model misadjustment (MMA) defined as $M_i:=\left\| \mathbf{a}_i-\mathbf{b%
}\right\| ^2$ is used as an objective convergence measure. Figure \ref
{figadapt} shows the MMA for different chaotic drive sequences obtained from
the LM as well as white Gausian noise. It is known that the best performance
is achieved for white Gausian noise. In this case the autocorrelation matrix
is $\mathbf{R}=\left( n+2\right) R_0^2\mathbf{I}$, with $R_0=E\left[
x^2\right] $. Comparing this with (\ref{corr1}) suggests that the zero-mean
chaotic sequence $\{x_i-\overline{x}\}$ generated by $f_4$ may achieve a
similar performance. This is confirmed by results in figure \ref{figadapt}.
In this simulation the damping parameter in (\ref{LMS2}) was set to the
maximum, i.e. $\mu =\frac{16}{3+2\,m}$ for the chaotic sequence and $\mu
=1/\left\| \mathbf{x}_i\right\| ^2$ for white noise. The results in figure 
\ref{figadapt} further confirm that the FIR DAF exhibits non-uniform
convergence when the chaotic sequence $\{x_i\}$ generated by $f_4$ is
applied. Moreover the MMA converges to a lower bound, which represents an
apparent performance limit. Whereas it is shown analytically in section \ref
{secLM} that the sequence $\{x_i-\overline{x}\}$ is uncorrelated, numerical
simulations show that the samples become correlated for $\lambda <4$, which
has consequences for the FIR DAF performance. This is revealed by the MMA
evolution for a sequence $\{x_i-\overline{x}\}$ generated by the LM $%
f_{3.95} $ shown in figure \ref{figadapt}. The effect of different $\lambda $
values is also visible in figure \ref{figadapt} for the test in which $%
\lambda $ is switched back and forth from 4 to 3.95. Obviously the
performance degrades when switched to $\lambda =3.95$ and becomes maximal
when switched back to $\lambda =4$. This observation is relevant for chaotic
coding schemes where $\lambda $ is varied according to the actual
information signal.

\begin{figure}[b]
\centerline{\psfig{file=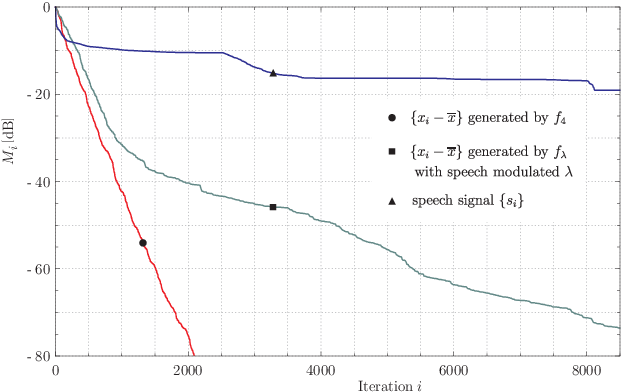,width=13.0cm}}
\caption{Adaptation performance when the FIR DAF is driven by a sequence $%
\left\{ x_i-\overline{x}\right\} $ generated by the LM $f_\lambda $ with
speech modulated bifurcation parameter according to $\lambda
_i=3.95+0.05\,s_i$, where $s_i$ is the speech signal in figure \ref
{figspeech}}
\label{adaptspeech}
\end{figure}
\begin{figure}[bh]
\centerline{\psfig{file=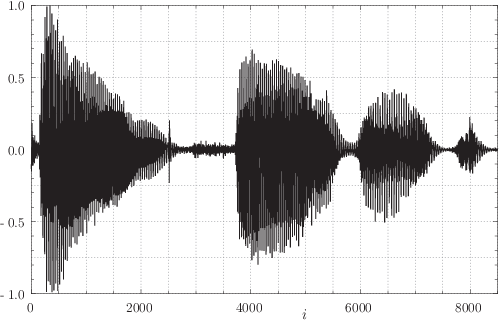,width=12.5cm}}
\caption{Normalized samples $s_i$ of speech signal sampled with 8 kHz
sampling rate.}
\label{figspeech}
\end{figure}

To analyze this effect the update performance is simulated for a modulated
bifurcation parameter according to $\lambda _i=3.95+0.05\,s_i$, where $s_i$
is the normalized speech signal in figure \ref{figspeech}, i.e. $|s|\leq 1$,
sampled with 8 kHz.

The resulting MMA is shown in figure \ref{adaptspeech}. Additionally shown
is the FIR DAF performance when the adaptation is driven by the original
speech signal $s_i$, which represents the classical approach \cite
{messer/ieee/84},\cite{haykin2001}. Although the chaotically coded speech
signal does not even lead close to the performance of white Gausian noise,
it leads to an adaptation that is significantly better than that achieved by
the original speech signal. The simulations results suggest that, when the
chaotic sequences generated by the LM are intended to be used as drive
signals for LMS FIR DAF channel estimation schemes, the chaotic modulation
schemes shall ensure that $\lambda $ remains close to 4.

\begin{figure}[tbh]
\centerline{\psfig{file=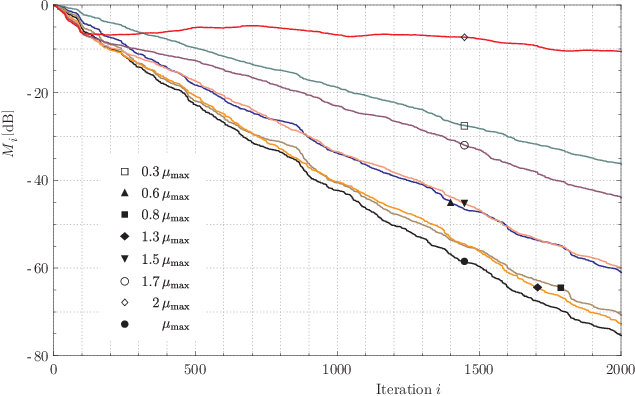,width=13.0cm}}
\caption{Adaptation performance of FIR DAF when driven by $\{x_i-\overline{x}%
\}$ generated by $f_4$ for different values of the damping parameter $\mu $
in the LMS scheme. Any deviation from the maximal value $\mu _{\max
}=16/\left( 3+2\,m\right) $ clearly impairs the performance.}
\label{figmu}
\end{figure}

In the preceding simulations the damping parameter has been set to the upper
bound estimated by (\ref{mumax}), and it is instructive to study the effect
of variations from this bound. The simulation results shown in figure \ref
{figmu} confirm that any variation from the upper bound (\ref{mumax})
impairs the performance. Moreover increasing $\mu $ eventually leads to
instability. It should be noticed that even for the upper bound the FIR DAF
may not be stable if the transmitted sequence is contaminated with additive
disturbing signals and noise, so that $\mu $ should be lower than the
maximum in (\ref{mumax}).

\vspace{3ex}%

\section{Summary}

The logistic map has been proposed for use in chaotic coding and spread
spectrum transmission systems. The statistical properties of sequences $%
\{x_i\}$ generated by the logistic map $f_4$ are deemed ideal for channel
estimation. This assumption was based solely on numerical simulation
results. In this contribution the higher order statistical moments and the
autocorrelation of the ergodic chaotic logistic map $f_4$ are derived. It is
proven that samples of the zero-mean sequence $\{x_i-\overline{x}\}$ are
uncorrelated and exhibit a flat power spectral density.

Upon these results the performance of FIR digital adaptive filters (DAF) is
analyzed when updated by a least mean squares (LMS) algorithm. It is shown
analytically that using zero-mean sequences of $f_4$ leads to the maximal
possible FIR DAF performance, which is for instance achieved by white
Gausian noise. An optimal value for the damping parameter in the LMS scheme
when driven by $\{x_i-\overline{x}\}$ is derived. These considerations are
confirmed by the presented simulation results.

It can be concluded that the ergodic logistic map $f_4$ dos in fact generate
sequences that are ideal for channel estimation and spread spectrum
applications. Since this does not strictly apply to the logistic map $%
f_\lambda $ with $\lambda <4$ it can also be concluded that for such
applications the chaotic modulation scheme must ensure that $\lambda $
remains close to 4. The analytical treatment of the statistics of the
chaotic $f_\lambda $ for $\lambda <4$ is still an open problem, however.

\end{document}